\documentclass[10pt]{article}
\usepackage{amssymb,amsmath}
\usepackage{graphicx}
\newtheorem{theorem}{Theorem}
\newtheorem{lemma}[theorem]{Lemma}

\begin{document}
\title{Note on generating all subsets of a finite set with disjoint unions}
\author{David Ellis}
\date{November 2008}
\maketitle

\begin{abstract}
We call a family \(\mathcal{G} \subset \mathbb{P}[n]\) a \(k\)-\textit{generator} of \(\mathbb{P}[n]\) if every \(x \subset [n]\) can be expressed as a union of at most \(k\) disjoint sets in \(\mathcal{G}\). Frein, L\'ev\^eque and Seb\H o \cite{Leveque} conjectured that for any \(n \geq k\), such a family must be at least as large as the \(k\)-generator obtained by taking a partition of \([n]\) into classes of sizes as equal as possible, and taking the union of the power-sets of the classes. We generalize a theorem of Alon and Frankl \cite{alon} in order to show that for fixed \(k\), any \(k\)-generator of \(\mathbb{P}[n]\) must have size at least \(k2^{n/k}(1-o(1))\), thereby verifying the conjecture asymptotically for multiples of \(k\).
\end{abstract}

\section{Introduction}
We call a family \(\mathcal{G} \subset \mathbb{P}[n]\) a \(k\)-\textit{generator} of \(\mathbb{P}[n]\) if every \(x \subset [n]\) can be expressed as a union of at most \(k\) \textrm{disjoint} sets in \(\mathcal{G}\). Frein, L\'ev\^eque and Seb\H o \cite{Leveque} conjectured that for any \(n \geq k\), such a family must be at least as large as the \(k\)-generator
\[\mathcal{F}_{n,k} := \bigcup_{i=1}^{k} \mathbb{P}V_{i} \setminus \{\emptyset\}\]
where \((V_{i})\) is a partition of \([n]\) into \(k\) classes of sizes as equal as possible. For \(k=2\), removing the disjointness condition yields the stronger conjecture of Erd\H os -- namely, if \(\mathcal{G} \subset \mathbb{P}[n]\) is a family such that any subset of \([n]\) is a union (not necessarily disjoint) of at most two sets in \(\mathcal{G}\), then \(\mathcal{G}\) is at least as large as
\[\mathcal{F}_{n,2} = \mathbb{P}V_{1} \cup \mathbb{P}V_{2} \setminus \{\emptyset\}\]
where \((V_{1},V_{2})\) is a partition of \([n]\) into two classes of sizes \(\lfloor n /2 \rfloor\) and \(\lceil n/2 \rceil\). We refer the reader to for example Furedi and Katona \cite{furedi} for some results around the Erd\H os conjecture. In fact, Frein, L\'ev\^eque and Seb\H o \cite{Leveque} made the analagous conjecture for all \(k\). (We call a family \(\mathcal{G} \subset \mathbb{P}[n]\) a \(k\)-\textit{base} of \(\mathbb{P}[n]\) if every \(x \subset [n]\) can be expressed as a union of at most \(k\) sets in \(\mathcal{G}\); they conjectured that for any \(k \leq n\), any \(k\)-base of \(\mathbb{P}[n]\) is at least as large as \(\mathcal{F}_{n,k}\).)

In this paper, we show that for \(k\) fixed, a \(k\)-generator must have size at least \(k2^{n/k}(1-o(1))\); when \(n\) is a multiple of \(k\), this is asymptotic to \(f(n,k) = |\mathcal{F}_{n,k}| = k(2^{n/k}-1)\). Our main tool is a generalization of a theorem of Alon and Frankl, proved via an Erdos-Stone type result.

We first remark that for a \(k\)-generator \(\mathcal{G}\), we have the following trivial bound on \(|\mathcal{G}|=m\). The number of ways of choosing at most \(k\) sets in \(\mathcal{G}\) must be at least the number of subsets of \([n]\), i.e.:
\[\sum_{i=0}^{k} {m \choose i} \geq 2^{n}\]
For fixed \(k\), the number of subsets of \([n]\) of size at most \(k-1\) is \(\sum_{i=0}^{k-1} {m \choose i} = \Theta(1/m) {m \choose k}\), so
\[\sum_{i=0}^{k} {m \choose i} = (1+\Theta(1/m)) {m \choose k} = (1+\Theta(1/m))m^{k}/k!\]
Hence,
\[m \geq (k!)^{1/k}2^{n/k}(1-o(1))\]
We will improve the constant from \((k!)^{1/k} \approx k/e\) to \(k\) by showing that for any fixed \(k \in \mathbb{N}\) and \(\delta > 0\), if \(m \geq 2^{(1/(k+1)+\delta)n}\), then any family \(\mathcal{G} \subset \mathbb{P}[n]\) of size \(m\) contains at most
\[\left(\frac{k!}{k^{k}}+o(1)\right){m \choose k}\]
unordered \(k\)-tuples \(\{A_{1},\ldots,A_{k}\}\) of pairwise disjoint sets, where the \(o(1)\) term tends to \(0\) as \(m \to \infty\) for fixed \(k,\delta\). In other words, if we consider the `Kneser graph' on \(\mathbb{P}[n]\), with edge set consisting of the disjoint pairs of subsets, the density of \(K_{k}\)'s in any sufficiently large \(\mathcal{G} \subset \mathbb{P}[n]\) is at most \(k!/k^{k} + o(1)\). (This generalizes Theorem 1.3 in \cite{alon}.) From the trivial bound above, any \(k\)-generator \(\mathcal{G} \subset \mathbb{P}[n]\) has size \(m \geq 2^{n/k}\), so putting \(\delta = 1/k(k+1)\), we will see that the number of unordered \(k\)-tuples of pairwise disjoint sets in \(\mathcal{G}\) is at most
\[\left(\frac{k!}{k^{k}}+o(1)\right){m \choose k}\]
so
\[2^{n} \leq \left(\frac{k!}{k^{k}}+o(1)+\Theta(1/m)\right){m \choose k} = \left(\frac{m}{k}\right)^{k}(1+o(1))\]
and therefore
\[m \geq k2^{n/k}(1-o(1))\]
where the \(o(1)\) term tends to 0 as \(n \to \infty\) for fixed \(k \in \mathbb{N}\).

\section{A preliminary Erd\H os-Stone type result}
We will need the following generalization of the Erd\H os-Stone theorem:
\begin{theorem}
Given \(r \leq s \in \mathbb{N}\) and \(\epsilon > 0\), if \(n\) is sufficiently large depending on \(r,s\) and \(\epsilon\), then any graph \(G\) on \(n\) vertices with at least \[\left(\frac{s(s-1)(s-2)\ldots(s-r+1)}{s^{r}}+\epsilon\right){n \choose r}\]
\(K_{r}\)'s contains a copy of \(K_{s+1}(t)\), where \(t \geq C_{r,s,\epsilon} \log n\) for some constant \(C_{r,s,\epsilon}\) depending on \(r,s,\epsilon\).
\end{theorem}
Note that the density \(\eta = \eta_{r,s} := \frac{s(s-1)(s-2)\ldots(s-r+1)}{s^{r}}\) above is the density of \(K_{r}\)'s in the \(s\)-partite Tur\'an graph with classes of size \(T\), \(K_{s}(T)\), when \(T\) is large.\\
\\
\textit{Proof:}\\
Let \(G\) be a graph with \(K_{r}\) density at least \(\eta + \epsilon\); let \(N\) be the number of \(l\)-subsets \(U \subset \mathcal{G}\) such that \(G[U]\) has \(K_{r}\)-density at least \(\eta + \epsilon/2\). Then, double counting the number of times an \(l\)-subset contains a \(K_{r}\),
\[N {l \choose r} + \left({n \choose r}-N\right)(\eta + \epsilon/2){l \choose r} \geq (\eta + \epsilon) {n \choose r}{n - r \choose l-r}\]
so rearranging,
\[N \geq \frac{\epsilon/2}{1-\eta-\epsilon/2} {n \choose l} \geq \tfrac{\epsilon}{2}{n \choose l}\]
Hence, there are at least \(\tfrac{\epsilon}{2}{n \choose l}\) \(l\)-sets \(U\) such that \(G[U]\) has \(K_{r}\)-density at least \(\eta + \epsilon/2\). But Erd\H os proved that the number of \(K_{r}\)'s in a \(K_{s+1}\)-free graph on \(l\) vertices is maximized by the \(s\)-partite Tur\'an graph on \(l\) vertices (Theorem 3 in \cite{erdos}), so provided \(l\) is chosen sufficiently large, each such \(G[U]\) contains a \(K_{s+1}\). Each \(K_{s+1}\) in \(G\) is contained in \({n-s-1 \choose l-s-1}\) \(l\)-sets, and therefore \(G\) contains at least
\[\tfrac{\epsilon}{2}\frac{{n \choose l}}{{ n - s-1 \choose l-s-1}} \geq \frac{\epsilon}{2} (n/l)^{s+1}\]
\(K_{s+1}\)'s, i.e. a positive density of \(K_{s+1}\)'s. Let \(a=s+1,\ c = \tfrac{\epsilon}{2l^{s+1}}\) and apply the following `blow up' theorem of Nikiforov (a slight weakening of Theorem 1 in \cite{nikiforov}):

\begin{theorem}
Let \(a \geq 2\), \(c^{a} \log n \geq 1\). Then any graph on \(n\) vertices with at least \(cn^{a}\) \(K_{a}\)'s contains a \(K_{a}(t)\) with \(t=\lfloor c^{a} \log n \rfloor\).
\end{theorem}

We see that provided \(n\) is sufficiently large depending on \(r,s\) and \(\epsilon\), \(G\) must contain a \(K_{s+1}(t)\) for \(t = \lfloor c^{s+1} \log n \rfloor = \lfloor (\tfrac{\epsilon}{2l^{s+1}})^{s+1} \log n \rfloor \geq C_{r,s,\epsilon} \log n\), proving Theorem 1. \(\square\)

\section{Density of \(K_{k}\)'s in large subsets of the Kneser graph}
We are now ready for our main result, a generalization of Theorem 1.3 in \cite{alon}:
\begin{theorem}
For any fixed \(k \in \mathbb{N}\) and \(\delta > 0\), if \(m \geq 2^{\left(\tfrac{1}{k+1}+\delta\right)n}\), then any family \(\mathcal{G} \subset \mathbb{P}[n]\) of size \(|\mathcal{G}| = m\) contains at most
\[\left(\frac{k!}{k^{k}}+o(1)\right){m \choose k}\]
unordered \(k\)-tuples \(\{A_{1},\ldots,A_{k}\}\) of pairwise disjoint sets, where the \(o(1)\) term tends to \(0\) as \(m \to \infty\) for fixed \(k,\delta\).
\end{theorem}
\textit{Proof:}\\
By increasing \(\delta\) if necessary, we may assume \(m=2^{\left(\tfrac{1}{k+1}+\delta\right)n}\). Consider the subgraph \(G\) of the `Kneser graph' on \(\mathbb{P}[n]\) induced on the set \(\mathcal{G}\), i.e. the graph \(G\) with vertex set \(\mathcal{G}\) and edge set \(\{xy: x \cap y = \emptyset\}\). Let \(\epsilon > 0\); we will show that if \(n\) is sufficiently large depending on \(k,\delta\) and \(\epsilon\), the density of \(K_{k}\)'s in \(G\) is less than \(\frac{k!}{k^{k}}+\epsilon\). Suppose the density of \(K_{k}\)'s in \(G\) is at least \(\frac{k!}{k^{k}}+\epsilon\); we will obtain a contradiction for \(n\) sufficiently large. Let \(l = m^{f}\) (we will choose \(f < \frac{\delta}{2(1+(k+1)\delta)}\) maximal such that \(m^{f}\) is an integer). By the argument above, there are at least \(\tfrac{\epsilon}{2}{m \choose l}\) \(l\)-sets \(U\) such that \(G[U]\) has \(K_{k}\)-density at least \(\frac{k!}{k^{k}} + \frac{\epsilon}{2}\). Provided \(m\) is sufficiently large depending on \(k,\delta\) and \(\epsilon\), by Theorem 1, each such \(G[U]\) contains a copy of \(K:=K_{k+1}(t)\) where \(t \geq C_{k,k,\epsilon/2} \log l = fC'_{k,\epsilon} \log m = C''_{k,\delta,\epsilon} \log m\). Any copy of \(K\) is contained in \({ m - (k+1)t \choose l-(k+1)t}\) \(l\)-sets, so \(G\) must contain at least \(\tfrac{\epsilon}{2}\frac{{m \choose l}}{{ m - (k+1)t \choose l-(k+1)t}} \geq \frac{\epsilon}{2} (m/l)^{(k+1)t}\) copies of \(K\).

But we also have the following lemma of Alon and Frankl (Lemma 4.3 in \cite{alon}), whose proof we include for completeness:

\begin{lemma}
\(G\) contains at most \((k+1)2^{n(1-\delta t)} {m \choose t}^{k+1} \frac{1}{(k+1)!}\) copies of \(K_{k+1}(t)\).
\end{lemma}
\textit{Proof:}\\
The probability that a \(t\)-subset \(\{A_{1},\ldots,A_{t}\}\) chosen uniformly at random from \(\mathcal{G}\) has union of size at most \(\tfrac{n}{k+1}\) is at most
\[\sum_{S \subset [n]: |S| \leq n/(k+1)} {2^{|S|} \choose t}/{m \choose t} \leq 2^{n} (2^{n/(k+1)}/m)^{t} = 2^{n(1-\delta t)}\]
Choose at random \(k+1\) such \(t\)-sets; the probability that at least one has union of size at most \(n/(k+1)\) is at most
\[(k+1)2^{n(1-\delta)t}\]
But this condition holds if our \(k+1\) \(t\)-sets are the vertex classes of a \(K_{k+1}(t)\) in \(G\). Hence, the number of copies of \(K_{k+1}(t)\) in \(G\) is at most
\[(k+1)2^{n(1-\delta t)} {m \choose t}^{k+1} \frac{1}{(k+1)!}\]
as required. \(\square\)\\

If \(m\) is sufficiently large depending on \(k,\delta\) and \(\epsilon\), we may certainly choose \(t \geq \lceil 4/\delta \rceil\), and comparing our two bounds gives
\[\tfrac{\epsilon}{2} (m/l)^{(k+1)t} \leq (k+1)2^{n(1-\delta t)}{m \choose t}^{k+1} \frac{1}{(k+1)!} \leq \tfrac{1}{2} 2^{n(1-\delta t)} m^{(k+1)t}\]
Substituting in \(l = m^{f}\), we get
\[\epsilon \leq 2^{n(1-\delta t)}m^{f(k+1)t}\]
Substituting in \(m = 2^{\left(\tfrac{1}{k+1}+\delta\right)n}\), we get
\[\epsilon \leq 2^{n(1-t(\delta - f(1+(k+1)\delta)))} \leq 2^{-n}\]
since we chose \(f < \frac{\delta}{2(1+(k+1)\delta)}\) and \(t \geq 4/\delta\). This is a contradiction if \(n\) is sufficiently large, proving Theorem 3. \(\square\)
\\

As explained above, our result on \(k\)-generators quickly follows:
\begin{theorem}
For fixed \(k \in \mathbb{N}\), any \(k\)-generator \(\mathcal{G}\) of \(\mathbb{P}[n]\) must contain at least \(k2^{n/k}(1-o(1))\) sets.
\end{theorem}
\textit{Proof:}\\
Let \(\mathcal{G}\) be a \(k\)-generator of \(\mathbb{P}[n]\), with \(|\mathcal{G}| = m\). As observed in the introduction, the trivial bound gives \(m \geq 2^{n/k}\), so applying Theorem 4 with \(\delta = 1/k(k+1)\), we see that the number of ways of choosing \(k\) pairwise disjoint sets in \(\mathcal{G}\) is at most
\[\left(\frac{k!}{k^{k}}+o(1)\right){m \choose k}\]
The number of ways of choosing less than \(k\) pairwise disjoint sets is, very crudely, at most \(\sum_{i=0}^{k-1} {m \choose i} = \Theta(1/m) {m \choose k}\); since every subset of \([n]\) is a disjoint union of at most \(k\) sets in \(\mathcal{G}\), we obtain
\[2^{n} \leq \left(\frac{k!}{k^{k}}+o(1)+\Theta(1/m)\right){m \choose k} = \left(\frac{m}{k}\right)^{k}(1+o(1))\]
(where the o(1) term tends to 0 as \(m \to \infty\)), and therefore
\[m \geq k2^{n/k}(1-o(1))\]
(where the \(o(1)\) term tends to 0 as \(n \to \infty\)). \(\square\)
\\
\\
\emph{Note:} The author wishes to thank Peter Keevash for bringing to his attention the result of Erd\H os in \cite{erdos},  after reading a previous draft of this paper in which a weaker, asymptotic version of Erd\H os' result was proved.

\end{document}